\documentclass[preprint,12pt]{elsarticle}
\usepackage{hyperref,url}
\usepackage{amsmath,amsthm,amssymb,bbm}

\newtheorem{theorem}{Theorem}[section]
\newtheorem{lemma}[theorem]{Lemma}
\newtheorem{corollary}[theorem]{Corollary}
\theoremstyle{definition}
\newtheorem{remark}[theorem]{Remark}
\newtheorem{definition}[theorem]{Definition}
\newtheorem{question}[theorem]{Question}
\begin{document}
\begin{frontmatter}
\title{Lower estimates for the norm and
the Kuratowski measure of noncompactness
of Wiener-Hopf type operators}

\author[OK]{Oleksiy Karlovych\corref{Oleksiy}}
\ead{oyk@fct.unl.pt}

\author[ES]{Eugene Shargorodsky}
\ead{eugene.shargorodsky@kcl.ac.uk}

\cortext[Oleksiy]{Corresponding author}

\address[OK]{%%
Centro de Matem\'atica e Aplica\c{c}\~oes,
Departamento de Matem\'atica,
Faculdade de Ci\^encias e Tecnologia,\\
Universidade Nova de Lisboa,
Quinta da Torre,
2829--516 Caparica,
Portugal}

\address[ES]{%
Department of Mathematics,
King's College London,
Strand, London WC2R 2LS,
United Kingdom}
%%%----------------------------------------------------------------------------
\begin{abstract}
Let $X(\mathbb{R}^n)$ be a Banach function space and 
$\Omega\subseteq\mathbb{R}^n$ be a measurable set of positive measure.
For a Fourier multiplier
$a$ on $X(\mathbb{R}^n)$, consider the Wiener-Hopf type operator
$W_\Omega(a):=r_\Omega F^{-1}aF e_\Omega$, where $F^{\pm 1}$ are the Fourier
transforms, $r_\Omega$ is the operator of restriction from $\mathbb{R}^n$
to $\Omega$ and $e_\Omega$ is the operator of extension by zero from $\Omega$
to $\mathbb{R}^n$. Let $X_2(\Omega)$ be the closure of 
$L^2(\Omega)\cap X(\Omega)$ in $X(\Omega)$. We show that if $X(\Omega)$
satisfies the so-called weak doubling property, then 
\[
\|a\|_{L^\infty(\mathbb{R}^n)}
\le \|W_\Omega(a)\|_{\mathcal{B}(X_2(\Omega),X(\Omega))}.
\]
Further, we prove that if $X(\Omega)$ satisfies the so-called separated 
doubling property, then the Kuratowski measure of noncompactness of 
$W_\Omega(a)$ admits the following lower estimate:
\[
\frac{1}{2}\|a\|_{L^\infty(\mathbb{R}^n)}
\le \|W_\Omega(a)\|_{\mathcal{B}(X_2(\Omega),X(\Omega)),\kappa}.
\]
These results are specified to the case of variable Lebesgue spaces
$L^{p(\cdot)}(C,w)$ with Muckenhoupt type weights $w$ over open cones
$C\subseteq\mathbb{R}^n$ with the vertex at the origin.
\end{abstract}

\end{frontmatter}
%%%------------------------------------------------------------------------
\section{Introduction}
For Banach spaces $\mathcal{X},\mathcal{Y}$, let 
$\mathcal{B}(\mathcal{X},\mathcal{Y})$ and 
$\mathcal{K}(\mathcal{X},\mathcal{Y})$ denote the Banach space of all
bounded linear operators from $\mathcal{X}$ to $\mathcal{Y}$
and its closed subspace consisting of all compact linear operators, 
respectively. If $\mathcal{X}=\mathcal{Y}$, then we will use the standard
abbreviations $\mathcal{B}(\mathcal{X}):=\mathcal{B}(\mathcal{X},\mathcal{X})$
and $\mathcal{K}(\mathcal{X}):=\mathcal{K}(\mathcal{X},\mathcal{X})$.

For a bounded subset $G$ of a Banach space $\mathcal{Y}$, the Kuratowski 
measure of noncompactness $\kappa(G)$ is defined as the greatest lower 
bound of the set of numbers $r$ such that $G$ can be covered by a finite 
family of sets of diameter at most $r$ 
(see, e.g., \cite[Definition~1.1.1]{AKPRS92} or \cite[Definition~5.1]{BM14}). 
It is clear that the equality $\kappa(G)=0$ implies precompactness
of $G\subset\mathcal{Y}$ (cf. \cite[Appendix~A4]{R91}).
For $A \in \mathcal{B}(\mathcal{X},\mathcal{Y})$, the Kuratowski
measure of noncompactness of $A$ is defined by
%%%
\begin{equation}\label{eq:Kuratowski}
\|A\|_{\mathcal{B}(\mathcal{X},\mathcal{Y}),\kappa} 
:= 
\frac12\,\kappa\left(A(B_\mathcal{X})\right) ,
\end{equation}
%%%
where $B_\mathcal{X}$ denotes the closed unit ball in $\mathcal{X}$. If 
the space
$\mathcal{X}$ is infinite-dimensional, then $\kappa(B_\mathcal{X})=2$
(see, e.g., \cite[Theorem~1.1.6]{AKPRS92}). This fact motivates
the choice of the factor $1/2$ in the above definition.

Similarly, the Hausdorff measure of noncompactness of a bounded
subset $G$ of a Banach space $\mathcal{Y}$ is the greatest lower bound
of the set of numbers $r>0$ such that $G$ can be covered by a finite family of
open balls of radius $r$ (see, e.g., \cite[Definition~1.1.2]{AKPRS92}
or \cite[Definition~5.2]{BM14}). The Hausdorff measure of noncompactness 
of $A\in\mathcal{B}(\mathcal{X},\mathcal{Y})$ is defined by
\[
\|A\|_{\mathcal{B}(\mathcal{X},\mathcal{Y}),\chi} :=\chi(A(B_{\mathcal{X}})).
\]
In view of \cite[Theorem~1.1.7] {AKPRS92}, 
\cite[Theorem~5.13 and Exercise~9 to Ch.~5]{BM14},
and \cite[inequality (3.29)]{LS71}, we have
\[
0\le 
\|A\|_{\mathcal{B}(\mathcal{X},\mathcal{Y}),\kappa}
\le
\|A\|_{\mathcal{B}(\mathcal{X},\mathcal{Y}),\chi}
\le 
\|A\|_{\mathcal{B}(\mathcal{X},\mathcal{Y}),\mathrm{e}},
\]
where 
\[
\|A\|_{\mathcal{B}(\mathcal{X},\mathcal{Y}),\mathrm{e}}
:=
\inf_{K\in\mathcal{K}(\mathcal{X},\mathcal{Y})} 
\|A-K\|_{\mathcal{B}(\mathcal{X},\mathcal{Y})},
\]
is the essential norm of $A$. It follows from \cite[Theorems~5.13 and~5.29]{BM14}
that
\begin{align*}
K\in\mathcal{K}(\mathcal{X},\mathcal{Y})
&\Longleftrightarrow
\|K\|_{\mathcal{B}(\mathcal{X},\mathcal{Y}),\kappa}=0
\\
&\Longleftrightarrow
\|K\|_{\mathcal{B}(\mathcal{X},\mathcal{Y}),\chi}=0
\\
&\Longleftrightarrow
\|K\|_{\mathcal{B}(\mathcal{X},\mathcal{Y}),\mathrm{e}}=0.
\end{align*}

Let $S(\mathbb{R}^n)$ denote the Schwartz space of rapidly
decaying infinitely differentiable functions, let
\[
(Fu)(\xi):=\widehat{u}(\xi):=\int_{\mathbb{R}^n}{u(x)}e^{-ix\xi}\,dx,
\quad
\xi\in\mathbb{R}^n,
\]
be the Fourier transform of $u\in S(\mathbb{R}^n)$, and let $F^{-1}$
denote the inverse Fourier transform on $S(\mathbb{R})$. The operators 
$F$ and $F^{-1}$ can be extended by continuity to $L^2(\mathbb{R}^n)$, 
and we will use the same symbols for their extensions to 
$L^2(\mathbb{R}^n)$.

Let $X(\mathbb{R}^n)$ be a Banach function space (see \cite{L55}
and also \cite[Ch.~1]{BS88}) and $X'(\mathbb{R}^n)$ be its associate 
space. We postpone the definitions of these spaces to 
Section~\ref{sec:BFS}. Here we only mention that the class of
Banach function spaces is very broad. It includes all Lebesgue 
spaces $L^p(\mathbb{R}^n)$ with $1\le p\le\infty$, all Orlicz spaces 
$L^\Phi(\mathbb{R}^n)$ generated by Young's functions $\Phi$
(see, e.g., \cite[Ch.~4, Section~8]{BS88}), all Lorentz spaces
$L^{p,q}(\mathbb{R}^n)$ with $1<p<\infty$ and $1\le q\le\infty$
(see, e.g., \cite[Ch.~4, Section~4]{BS88}), 
all Morrey spaces $\mathcal{M}_q^p(\mathbb{R}^n)$ with
$1\le q\le p<\infty$ (see \cite[Section~1.2.6]{SDH20}).
All the above spaces are translation-invariant. The class of Banach
function spaces also includes all variable Lebesgue 
spaces $L^{p(\cdot)}(\mathbb{R}^n)$ (see \cite{CF13,DHHR11}),
as well as weighted analogues of the above-mentioned spaces
(under some conditions on weights, see Lemma~\ref{le:wBFS} below).
These spaces are not translation-invariant.

A function $a\in L^\infty(\mathbb{R}^n)$ is said to belong to the set
$\mathcal{M}^0_{X(\mathbb{R}^n)}$ of Fourier multipliers on 
a Banach function space $X(\mathbb{R}^n)$ if
\[
\|a\|_{\mathcal{M}^0_{X(\mathbb{R}^n)}}
:=
\sup\left\{
\frac{\|F^{-1} aFu\|_{X(\mathbb{R}^n)}}{\|u\|_{X(\mathbb{R}^n)}}:
u \in \left(L^2(\mathbb{R}^n)\cap X(\mathbb{R}^n)\right) \setminus\{0\}
\right\}<\infty.
\]

Let $X(\mathbb{R}^n)$ be a Banach function space and 
$\Omega\subseteq\mathbb{R}^n$ be a measurable set of positive measure. 
Denote by $r_\Omega$ the operator of restriction from $\mathbb{R}^n$ 
to $\Omega$ and by $e_\Omega$ the operator of extension by $0$ 
from $\Omega$ to $\mathbb{R}^n$. Let $X(\Omega)$ be the space
of all measurable functions $f:\Omega\to\mathbb{C}$ such that
\[
\|f\|_{X(\Omega)}:=\|e_\Omega f\|_{X(\mathbb{R}^n)}<\infty.
\]
For $a \in \mathcal{M}^0_{X(\mathbb{R}^n)}$, let
\[
W_\Omega(a) u := r_\Omega F^{-1} aFe_\Omega u, 
\quad u \in L^2(\Omega)\cap X(\Omega) 
\]
be an operator of Wiener-Hopf type. Put
\[
\|W_\Omega(a)\|_{[X(\Omega)]} := 
\sup\left\{
\frac{\|W_\Omega(a) u\|_{X(\Omega)}}{\|u\|_{X(\Omega)}}:
u \in \left(L^2(\Omega)\cap X(\Omega)\right) \setminus\{0\}
\right\} .
\]
Let $X_2(\Omega)$ be the closure of $L^2(\Omega)\cap X(\Omega)$ in $X(\Omega)$.
Then the bounded linear operator 
\[
W_\Omega(a):L^2(\Omega)\cap X(\Omega)\to X(\Omega) 
\]
can be extended by continuity to a bounded linear operator from $X_2(\Omega)$ 
to $X(\Omega)$, which we will denote again by $W_\Omega(a)$. Moreover,
\[
\|W_\Omega(a)\|_{[X(\Omega)]}
=
\|W_\Omega(a)\|_{\mathcal{B}(X_2(\Omega),X(\Omega))}.
\]

We are interested in the relations between the operator norm,
the essential norm, and the Kuratowski and Hausdorff measures of noncompactness
of $W_{\mathbb{R}^n}(a)$ and $W_\Omega(a)$ with $\Omega\subset\mathbb{R}^n$
on fairly general Banach function spaces. Let us recall some results 
motivating our work. If $X(\mathbb{R}^n)$ is a separable translation-invariant
Banach function space and $a\in\mathcal{M}_{X(\mathbb{R}^n)}^0$,
then it follows from \cite[Theorem~5.1]{KS24} that
%%%
\begin{equation}\label{eq:multiplier-operator}
\|W_{\mathbb{R}^n}(a)\|_{\mathcal{B}(X(\mathbb{R}^n))}
=
\|W_{\mathbb{R}^n}(a)\|_{\mathcal{B}(X(\mathbb{R}^n)),\mathrm{e}}
=
\|W_{\mathbb{R}^n}(a)\|_{\mathcal{B}(X(\mathbb{R}^n)),\chi}.
\end{equation}
%%%
On the other hand, if $n=1$, $X(\mathbb{R})$ is a separable translation-invariant
Banach function space, and $a\in\mathcal{M}_{X(\mathbb{R})}^0$, then
it follows from \cite[Theorem~1.1]{KS-JMS} that
%%%
\begin{equation}\label{eq:WH-operator}
\|W_{\mathbb{R}_+}(a)\|_{\mathcal{B}(X(\mathbb{R}_+))}
=
\|W_{\mathbb{R}_+}(a)\|_{\mathcal{B}(X(\mathbb{R}_+)),\mathrm{e}}
=
\|W_{\mathbb{R}_+}(a)\|_{\mathcal{B}(X(\mathbb{R}_+)),\chi}.
\end{equation}
%%%
These equalities extended some earlier results presented in \cite[Section~9.5]{BS06}
and \cite[Theorem~1]{KV25}. They imply that if $n=1$ and $X(\mathbb{R})$
is a separable translation-invariant Banach function space, then 
%%%
\begin{equation}\label{eq:maximal-noncompactness-corollaries}
W_\mathbb{R}(a)\in\mathcal{K}(X(\mathbb{R}))
\ \Longleftrightarrow\
a=0,
\quad
W_{\mathbb{R}_+}(a)\in\mathcal{K}(X(\mathbb{R}_+))
\ \Longleftrightarrow\ 
a=0.
\end{equation}
%%%

Our main motivation for this work is to answer the following question: 
%%%----------------------------------------------------------------------------
\begin{question}\label{question:main}
Are the equivalences in \eqref{eq:maximal-noncompactness-corollaries} 
true in the setting of non-trans\-la\-tion-invariant Banach function 
spaces, like weighted Lebesgue spaces $L^p(\mathbb{R},w)$ or variable 
Lebesgue spaces $L^{p(\cdot)}(\mathbb{R})$ (see, e.g., \cite{CF13,DHHR11} 
and Section~\ref{sec:weighted-variable-Lebesgue-spaces})? 
\end{question}
%%%----------------------------------------------------------------------------
The answer to this question is important in various problems of the 
Fredholm theory of Fourier convolution operators and Wiener-Hopf operators.

We are not aware of quantitative estimates for the essential norm
or Kuratowski measure of noncompactness of operators $W_{\mathbb{R}}(a)$
or $W_{\mathbb{R}_+}(a)$ on non-translation-invariant spaces.
As for the qualitative results, the first equivalence 
in \eqref{eq:maximal-noncompactness-corollaries} was established in the case
of a separable Banach function space $X(\mathbb{R})$ such that the 
Hardy-Littlewood maximal operator $M$ is bounded on $X(\mathbb{R})$ and on 
its associate space $X'(\mathbb{R})$ (see \cite[Theorem~1]{FKK19}). 
Recall that the Hardy-Littewood maximal operator on $\mathbb{R}^n$
is defined by
\[
(Mf)(x):=\sup_{B\ni x}\frac{1}{|B|}\int_B |f(y)|\,dy,
\quad
x\in\mathbb{R}^n,
\]
where the supremum is taken over all balls in 
$B\subseteq\mathbb{R}^n$ containing $x$.

Our first result is the following complement to \cite[Theorem~1.3]{KS19}.
%%%----------------------------------------------------------------------------
\begin{theorem}\label{th:norm-estimate} 
Let $X(\mathbb{R}^n)$ be a Banach function space and 
$\Omega\subseteq\mathbb{R}^n$ be a measurable set of infinite measure
such that $X(\Omega)$ has the weak doubling property 
(see Definition~\ref{def:weak-doubling} below).
If $a \in \mathcal{M}^0_{X(\mathbb{R}^n)}$, then
\[
\|a\|_{L^\infty(\mathbb{R}^n)} 
\le 
\|W_\Omega(a)\|_{\mathcal{B}(X_2(\Omega),X(\Omega))}.
\]
\end{theorem}
%%%----------------------------------------------------------------------------
Our main result is the following lower estimate for the Kuratowski measure
of noncompactness of the Wiener-Hopf type operator $W_\Omega(a)$.
%%%----------------------------------------------------------------------------
\begin{theorem}\label{th:Kuratowski-estimate} 
Let $X(\mathbb{R}^n)$ be a Banach function space and 
$\Omega\subseteq\mathbb{R}^n$ be a measurable set of infinite measure
such that $X(\Omega)$ has the separated doubling property 
(see Definition~\ref{def:separated-doubling} below). 
If $a \in \mathcal{M}^0_{X(\mathbb{R}^n)}$, then
\[
\frac12\, \|a\|_{L^\infty(\mathbb{R}^n)} 
\le 
\|W_\Omega(a)\|_{\mathcal{B}(X_2(\Omega),X(\Omega)),\kappa} .
\]
\end{theorem}
%%%----------------------------------------------------------------------------

It would be interesting to know whether the constant $1/2$ in the above
inequality is sharp or it can be substituted by a larger constant 
$c\in(1/2,1]$.

Let us formulate a corollary of the above results in a form convenient for
applications.

An open subset $C\subseteq\mathbb{R}^n$ is said to be an open cone with 
the vertex at the origin if $x\in C$ implies that $\alpha x\in C$ for all $\alpha>0$.
%%%----------------------------------------------------------------------------
\begin{corollary}\label{co:nice-for-cones}
Let $X(\mathbb{R}^n)$ be a separable Banach function space with the associate 
space $X'(\mathbb{R}^n)$ and let $C\subseteq \mathbb{R}^n$ be an open cone 
with the vertex at the origin. Suppose that
%%%
\begin{equation}\label{eq:Berezhnoi}
\sup_{B}\frac{1}{|B|}
\|\chi_B\|_{X(\mathbb{R}^n)}
\|\chi_B\|_{X'(\mathbb{R}^n)}
<\infty,
\end{equation}
%%%
where the supremum is taken over all balls $B\subseteq\mathbb{R}^n$.
If $a \in \mathcal{M}^0_{X(\mathbb{R}^n)}$, then
\[
\|a\|_{L^\infty(\mathbb{R}^n)} 
\le 
\|W_C(a)\|_{\mathcal{B}(X(C))},
\quad
\frac12\, \|a\|_{L^\infty(\mathbb{R}^n)} 
\le 
\|W_C(a)\|_{\mathcal{B}(X(C)),\kappa} .
\]
In particular,
%%%
\begin{equation}\label{eq:nice-for-cones}
W_{\mathbb{R}^n}(a)\in\mathcal{K}(X(\mathbb{R}^n))
\ \Longleftrightarrow\
a=0,
\quad
W_{C}(a)\in\mathcal{K}(X(C))
\ \Longleftrightarrow\ 
a=0.
\end{equation}
\end{corollary}
%%%----------------------------------------------------------------------------

This corollary follows immediately from 
Theorems~\ref{th:norm-estimate}--\ref{th:Kuratowski-estimate}
and Lemmas~\ref{le:X2=X} and~\ref{le:doubling-properties-over-cones} below.

As far as we know, condition \eqref{eq:Berezhnoi} was introduced by Berezhnoi
\cite{B99}. It follows from \cite[Lemma~3.2]{H12} that if the Hardy-Littlewood
maximal operator $M$ is bounded on a Banach function space
$X(\mathbb{R}^n)$ or on its associate space $X'(\mathbb{R}^n)$, then
\eqref{eq:Berezhnoi} holds. We will explain in 
Section~\ref{sec:weighted-variable-Lebesgue-spaces}
that condition \eqref{eq:Berezhnoi} is strictly weaker than the boundedness
of $M$ on $X(\mathbb{R}^n)$ on the class of all separable Banach
function spaces $X(\mathbb{R}^n)$. So, the first equivalence in 
\eqref{eq:nice-for-cones} extends \cite[Theorem~1.1]{FKK19} for $n=1$.

The paper is organised as follows. In Section~\ref{sec:BFS}, we recall
the definition of a Banach function space $X(\mathbb{R}^n)$ and its associate
space $X'(\mathbb{R}^n)$. In Section~\ref{sec:doubling}, we discuss
the weak doubling and separated doubling properties and show that
if condition \eqref{eq:Berezhnoi} is satisfied then the space $X(C)$
has both properties for every cone $C$ with the vertex at the origin.
Theorems~\ref{th:norm-estimate} and~\ref{th:Kuratowski-estimate} are proved
in Sections~\ref{sec:norm-estimate-proof} 
and~\ref{sec:Kuratowski-estimate-proof}, respectively. 
Finally, in Section~\ref{sec:weighted-variable-Lebesgue-spaces},
we specify Corollary~\ref{co:nice-for-cones} to the case of weighted
variable Lebesgue spaces (including the case of Lebesgue spaces with 
Muckenhoupt weights).
%%%----------------------------------------------------------------------------
\section{Banach function spaces}\label{sec:BFS}
The set of all Lebesgue measurable complex-valued functions on $\mathbb{R}^n$ 
is denoted by $\mathfrak{M}(\mathbb{R}^n)$. Let $\mathfrak{M}^+(\mathbb{R}^n)$ 
be the subset of functions in $\mathfrak{M}(\mathbb{R}^n)$ whose
values lie  in $[0,\infty]$. For a measurable set $E\subseteq\mathbb{R}^n$, 
its Lebesgue measure and the characteristic function are denoted by $|E|$ and
$\chi_E$, respectively. Following \cite[p.~3]{L55} (see also 
\cite[Ch.~1, Definition~1.1]{BS88} and \cite[Definition~6.1.5]{PKJF13}), a 
mapping $\rho:\mathfrak{M}^+(\mathbb{R}^n)\to [0,\infty]$ is called a Banach 
function norm if, for all functions $f,g, f_n \ (n\in\mathbb{N})$ in 
$\mathfrak{M}^+(\mathbb{R}^n)$, for all constants $a\ge 0$, and for all 
measurable subsets $E$ of $\mathbb{R}^n$, the following properties hold:
\begin{eqnarray*}
{\rm (A1)} &\quad & \rho(f)=0  \Leftrightarrow  f=0\ \mbox{a.e.}, \
\rho(af)=a\rho(f), \
\rho(f+g) \le \rho(f)+\rho(g),\\
{\rm (A2)} &\quad &0\le g \le f \ \mbox{a.e.} \ \Rightarrow \ \rho(g)
\le \rho(f)
\quad\mbox{(the lattice property)},
\\
{\rm (A3)} &\quad &0\le f_n \uparrow f \ \mbox{a.e.} \ \Rightarrow \
       \rho(f_n) \uparrow \rho(f)\quad\mbox{(the Fatou property)},\\
{\rm (A4)} &\quad & 
{E \text{ is bounded}}  
\Rightarrow \rho(\chi_E) <\infty,\\
{\rm (A5)} &\quad & 
{E \text{ is bounded}}  
\Rightarrow \int_E f(x)\,dx \le C_E\rho(f),
\end{eqnarray*}
%%%%
where $C_E \in (0,\infty)$ may depend on $E$ and $\rho$ but is
independent of $f$.  When functions differing only on a set of measure zero
are identified, the set $X(\mathbb{R}^n)$ of functions 
$f\in\mathfrak{M}(\mathbb{R}^n)$ for which $\rho(|f|)<\infty$ is called a 
Banach function space. For each $f\in X(\mathbb{R}^n)$, the norm of $f$ is 
defined by
\[
\left\|f\right\|_{X(\mathbb{R}^n)} :=\rho(|f|).
\]
With this norm and under natural linear space operations, the set 
$X(\mathbb{R}^n)$ becomes a Banach space (see 
\cite[Ch.~1, \S1, Theorem~1]{L55} or 
\cite[Ch.~1, Theorems~1.4 and~1.6]{BS88}). 
If $\rho$ is a Banach function norm, its associate norm $\rho'$ is defined on
$\mathfrak{M}^+(\mathbb{R}^n)$ by
\[
\rho'(g):=\sup\left\{
\int_{\mathbb{R}^n} f(x)g(x)\,dx \ : \ f\in \mathfrak{M}^+(\mathbb{R}^n), 
\ \rho(f) \le 1 \right\}, 
\quad 
g\in \mathfrak{M}^+(\mathbb{R}^n).
\]
Then $\rho'$ is itself a Banach function norm (see \cite[Ch.~1, \S1]{L55} 
or \cite[Ch.~1, Theorem~2.2]{BS88}). The Banach function space 
$X'(\mathbb{R}^n)$ determined by the Banach function norm $\rho'$ is called 
the associate space (K\"othe dual) of $X(\mathbb{R}^n)$. The associate space 
$X'(\mathbb{R}^n)$ can be identified with a subspace of the (Banach) dual
space $[X(\mathbb{R} ^n)]^*$.
%%%----------------------------------------------------------------------------
\begin{remark}
We note that our definition of a Banach function space is slightly 
different from that found in \cite[Chap.~1, Definition~1.1]{BS88}
and \cite[Definition~6.1.5]{PKJF13}. 
In particular, in Axioms (A4) and (A5)
we assume that the set $E$ is a bounded set, whereas it is sometimes 
assumed that $E$ merely satisfies $|E| <\infty$. We do this so that the 
weighted Lebesgue spaces with Muckenhoupt weights satisfy
Axioms (A4) and (A5). The same applies to Morrey spaces
$\mathcal{M}_q^p(\mathbb{R}^n)$ with $1\le q\le p<\infty$
(see \cite[Section~1.2.6]{SDH20}). Moreover, it is well known that
all main elements of the general theory of Banach function spaces
work with (A4) and (A5) stated for bounded sets \cite{L55} (see also 
the discussion at the beginning of Chapter~1 on page~2 of \cite{BS88} 
and \cite{LN24}).
Unfortunately, we overlooked that the definition of a Banach function 
space in our previous works \cite{FKK19,KS19,KS14}
had to be changed by replacing in Axioms (A4) and (A5)
the requirement of $|E|<\infty$ by the requirement 
that $E$ is a bounded set to include weighted Lebesgue spaces
in our considerations. However, the results proved in the above
papers remain true. 
\end{remark}
%%%----------------------------------------------------------------------------
Let $X(\mathbb{R}^n)$ be a Banach function space generated by a Banach function 
norm $\rho$. We say that $f\in X_{\rm loc}(\mathbb{R}^n)$ if 
$f\chi_E\in X(\mathbb{R}^n)$ for every bounded measurable set 
$E\subset\mathbb{R}^n$. Further, by $C_0^\infty(\mathbb{R}^n)$ we denote the 
set of all infinitely differentiable functions with compact support.
A function $w:\mathbb{R}^n\to[0,\infty]$ is said to 
be a weight if $0<w<\infty$ a.e. on $\mathbb{R}^n$.
The proof of the following lemma given in \cite[Lemmas~2.4, 2.9, and 2.12]{KS14}
for $n=1$ can be easily extended to arbitrary $n\ge 1$.
%%%----------------------------------------------------------------------------
\begin{lemma}
\label{le:wBFS}
Let $X(\mathbb{R}^n)$ be a Banach function space generated by a Banach function 
norm $\rho$, let $X'(\mathbb{R}^n)$ be its associate space, and let 
$w:\mathbb{R}\to[0,\infty]$ be a weight. Suppose that 
$w\in X_{\rm loc}(\mathbb{R}^n)$ and $w^{-1}\in X_{\rm loc}'(\mathbb{R}^n)$. 

\begin{enumerate}
\item[{\rm(a)}]
Then 
\[
\rho_w(f):=\rho(fw),\quad f\in\mathfrak{M}^+(\mathbb{R}^n),
\]
is a Banach function norm and
\[
X(\mathbb{R}^n,w):=\{f\in \mathfrak{M}(\mathbb{R}^n):fw\in X(\mathbb{R}^n)\}
\]
is a Banach function space generated by the Banach function norm $\rho_w$.
The space $X'(\mathbb{R}^n,w^{-1})$ is the associate space of $X(\mathbb{R}^n,w)$.

\item[{\rm (b)}]
If $X(\mathbb{R}^n)$ is separable (resp., reflexive), then 
$X(\mathbb{R}^n,w)$ is separable (resp., reflexive).

\item[{\rm (c)}]
If $X(\mathbb{R}^n)$ is separable, then the set $C_0^\infty(\mathbb{R}^n)$
is dense in the space $X(\mathbb{R}^n,w)$.
\end{enumerate}
\end{lemma}
%%%----------------------------------------------------------------------------
We conclude this section with a corollary of part (c) of the above lemma.
%%%----------------------------------------------------------------------------
\begin{lemma}\label{le:X2=X}
Let $X(\mathbb{R}^n)$ be a Banach function space and let 
$\Omega\subseteq\mathbb{R}^n$ be a measurable set of positive measure.
If $X(\mathbb{R}^n)$ is separable, then $X_2(\Omega)=X(\Omega)$.
\end{lemma}
%%%----------------------------------------------------------------------------
\begin{proof}
Let $f\in X(\Omega)$. Then 
$g:=e_\Omega f\in X(\mathbb{R}^n)$. It follows from Lemma~\ref{le:wBFS}(c)
that for every $\varepsilon>0$ there exists 
$h\in C_0^\infty(\mathbb{R}^n)\subseteq L^2(\mathbb{R}^n)\cap X(\mathbb{R}^n)$
such that $\|g-h\|_{X(\mathbb{R}^n)}<\varepsilon$. Then
$r:=r_\Omega h\in L^2(\Omega)\cap X(\Omega)$ and
%%%
\begin{align*}
\|f-r\|_{X(\Omega)}
&=
\|r_\Omega e_\Omega f- r_\Omega h\|_{X(\Omega)}
=
\|r_\Omega(g-h)\|_{X(\Omega)}
\\
&=
\|e_\Omega r_\Omega(g-h)\|_{X(\mathbb{R}^n)}
=
\|\chi_\Omega(g-h)\|_{X(\mathbb{R}^n)}
\le 
\|g-h\|_{X(\mathbb{R}^n)}
<
\varepsilon.
\end{align*}
%%%
Hence $f\in X_2(\Omega)$. So, $X(\Omega)\subseteq X_2(\Omega)$.
The reverse inclusion $X_2(\Omega)\subseteq X(\Omega)$ is obvious.
\end{proof}
%%%----------------------------------------------------------------------------
\section{Doubling properties}\label{sec:doubling}
Let $\Omega\subseteq\mathbb{R}^n$ be a measurable set of infinite measure
and $X(\mathbb{R}^n)$ be a Banach function space. For $y\in\mathbb{R}^n$ and 
$R>0$, let $B(y,R) := \{x \in \mathbb{R}^n : |x-y|< R\}$ be the open ball of 
radius $R$ centred at $y$. 

According to \cite[Definition~3.2]{KS19}, a Banach function space 
$X(\mathbb{R}^n)$ is said to have the strong doubling property if there 
exist a number $\tau > 1$ and a constant $A_\tau > 0$ such that for all 
$R > 0$ and $y \in \mathbb{R}^n$,
%%%
\begin{equation}\label{eq:strong-doubling}
\frac{\|\chi_{B(y,\tau R)}\|_{X(\mathbb{R}^n)}}
{\|\chi_{B(y, R)}\|_{X(\mathbb{R}^n)}} 
\le 
A_\tau.
\end{equation}

Now we introduce two doubling properties related to the space $X(\Omega)$. 
The first one is an adaptation of the weak doubling property for 
$X(\mathbb{R}^n)$ (see \cite[Definition~1.2]{KS19}).
%%%----------------------------------------------------------------------------
\begin{definition}\label{def:weak-doubling}
The space $X(\Omega)$ is said to have the weak doubling property if there 
exists a number $\tau > 1$  such that
%%%
\begin{equation}\label{eq:weak-doubling}
D_{X(\Omega), \tau} 
:= 
\liminf_{R \to \infty}
\left(\inf_{B(y,\tau R) \subseteq \Omega}
\frac{\|\chi_{B(y,\tau R)}\|_{X(\Omega)}}
{\|\chi_{B(y, R)}\|_{X(\Omega)}}
\right)<\infty.
\end{equation}
\end{definition}
%%%----------------------------------------------------------------------------
The strong doubling property is considerably stronger 
than the weak doubling property.
Indeed, it is easy to see that \eqref{eq:weak-doubling} is equivalent to the 
following: for every $\varepsilon > 0$, there exist $R_j \nearrow \infty$ 
and  $y_j \in \Omega$, $j \in \mathbb{N}$ such that 
$B(y_j,\tau R_j) \subseteq \Omega$ and
\[
\frac{\|\chi_{B(y_j,\tau R_j)}\|_{X(\Omega)}}
{\|\chi_{B(y_j, R_j)}\|_{X(\Omega)}} 
\le 
D_{X(\Omega), \tau} + \varepsilon , 
\quad 
j \in \mathbb{N} .
\]
So, the difference between the strong doubling property 
and the weak doubling property is that the former requires 
estimate \eqref{eq:strong-doubling} 
to hold for all balls, while the latter requires it 
to hold only for some sequence of balls 
with radii going to infinity.

Now we introduce the separated doubling property of $X(\Omega)$ 
playing a fundamental role in our main result 
(Theorem~\ref{th:Kuratowski-estimate}).
%%%----------------------------------------------------------------------------
\begin{definition}\label{def:separated-doubling}
The space $X(\Omega)$ is said to have the separated 
doubling property if there exist $\tau > 1$,  $C_\tau  \in (0, \infty)$,  
$R_j \nearrow \infty$ and $y_j \in \Omega$, $j \in \mathbb{N}$
such that 
\[
B(y_j,\tau R_j) \subseteq \Omega,
\quad  
B(y_j,\tau R_j) \cap B(y_k,\tau R_k) = \emptyset 
\quad\mbox{for}\quad k \not= j, 
\]
and
\[
\frac{\|\chi_{B(y_j,\tau R_j)}\|_{X(\Omega)}}
{\|\chi_{B(y_j, R_j)}\|_{X(\Omega)}} 
\le 
C_\tau , 
\quad 
j \in \mathbb{N} .
\]
We denote by $S_{X(\Omega), \tau}$ the infimum of the constants $C_\tau$ for 
which the above conditions are satisfied.
\end{definition}
%%%----------------------------------------------------------------------------
It is clear that if $X(\Omega)$ has the separated doubling property, then it 
also has the weak doubling property. It would be interesting to clarify 
whether there exist a Banach function space $X(\mathbb{R}^n)$ and
a measurable set $\Omega\subseteq\mathbb{R}^n$ such that $X(\Omega)$
has the weak doubling property but does not have the separated doubling 
property.

We will need the following analogue of \cite[Lemma~3.1]{KS19}.
%%%----------------------------------------------------------------------------
\begin{lemma}\label{le:=1} 
Let $X(\mathbb{R}^n)$ be a Banach function space and 
$\Omega\subseteq\mathbb{R}^n$ be a measurable set of infinite measure.
%%%
\begin{enumerate}
\item[{\rm (a)}]
If $X(\Omega)$ has the weak doubling property, then
\[
\inf_{\tau > 1} D_{X(\Omega), \tau} = 1.
\]

\item[{\rm(b)}]
If $X(\Omega)$ has the separated doubling property, then
\[
\inf_{\tau > 1} S_{X(\Omega), \tau} = 1.
\]
\end{enumerate}
\end{lemma}
%%%----------------------------------------------------------------------------
\begin{proof}
(a) Since the proof is almost identical to that of \cite[Lemma~3.1]{KS19},
it is omitted here.

(b) The proof is similar to that of \cite[Lemma 3.1]{KS19}.
Since $X(\Omega)$ satisfies the separated doubling property, there exists a
number $\varrho>1$ such that $S_{X(\Omega), \varrho}<\infty$. Then there exist
$R_j \nearrow \infty$ and $y_j \in \Omega$, $j \in \mathbb{N}$
such that $B(y_j, \varrho R_j) \subseteq \Omega$,  
$B(y_j,\varrho R_j) \cap B(y_k,\varrho R_k) = \emptyset$ for $k \not= j$, and
\[
\frac{\|\chi_{B(y_j,\varrho R_j)}\|_{X(\Omega)}}
{\|\chi_{B(y_j, R_j)}\|_{X(\Omega)}} \le S_{X(\Omega), \varrho} + 1 , 
\qquad 
j \in \mathbb{N} .
\]
Assume, contrary to the hypothesis, that
\[
S := \inf_{\tau > 1} S_{X(\Omega), \tau} > 1.
\]
Take an arbitrary $N \in \mathbb{N}$ and consider $\tau = \varrho^{1/N}$.
Then
%%%
\begin{equation}\label{eq:>1}
S_{X(\Omega), \tau} \ge S > S_0 := \frac{S + 1}{2} > 1.
\end{equation}
If there exists a strictly increasing sequence of natural numbers $j_k$, 
$k \in \mathbb{N}$, such that 
\[
\frac{\|\chi_{B(y_{j_k},\tau R_{j_k})}\|_{X(\Omega)}}
{\|\chi_{B(y_{j_k}, R_{j_k})}\|_{X(\Omega)}} < S_0 ,
\]
then
it follows from the definition of $S_{X(\Omega), \tau}$ that 
$S_{X(\Omega), \tau} \le S_0$, which contradicts \eqref{eq:>1}.
Hence, there exists $j_0 \in \mathbb{N}$ such that 
\[
\frac{\|\chi_{B(y_j,\tau R_j)}\|_{X(\Omega)}}
{\|\chi_{B(y_j, R_j)}\|_{X(\Omega)}} \ge S_0  \quad\mbox{for all}  \quad j \ge j_0 .
\]
Similarly, one can prove the existence of $J \in \mathbb{N}$ such that 
\[
\frac{\|\chi_{B(y_j,\tau^\ell R_j)}\|_{X(\Omega)}}
{\|\chi_{B(y_j, \tau^{\ell - 1} R_j)}\|_{X(\Omega)}} \ge S_0  
\quad
\mbox{for all}  \quad j \ge J , \ \ell = 1, \dots, N .
\]
Then
\[
\frac{\|\chi_{B(y_j,\varrho R_j)}\|_{X(\Omega)}}
{\|\chi_{B(y_j, R_j)}\|_{X(\Omega)}}
=
\prod_{\ell = 1}^N
\frac{\|\chi_{B(y_j,\tau^\ell R_j)}\|_{X(\Omega)}}
{\|\chi_{B(y_j, \tau^{\ell - 1} R_j)}\|_{X(\Omega)}}
\ge
S_0^N , \qquad j \ge J .
\]
Thus, $S_{X(\Omega), \varrho} + 1 \ge S_0^N$ for all $N \in \mathbb{N}$,
which is impossible since $S_0 > 1$ and $S_{X(\Omega), \varrho}  < \infty$.
The obtained contradiction completes the proof.
\end{proof}
%%%----------------------------------------------------------------------------
Now we give a simple sufficient condition guaranteeing that the space
$X(\Omega)$ has the separated doubling property and, hence, the weak
doubling property.
%%%----------------------------------------------------------------------------
\begin{lemma}\label{le:doubling-properties-over-cones}
Let $X(\mathbb{R}^n)$ be a Banach function space satisfying 
\eqref{eq:Berezhnoi}.
Then $X(C)$ has the weak doubling property and the separated doubling
property for every open cone $C$ with the vertex at the origin.
\end{lemma}
%%%----------------------------------------------------------------------------
\begin{proof}
It follows from \cite[Lemma~3.3]{KS19} that condition \eqref{eq:Berezhnoi}
implies that the space $X(\mathbb{R}^n)$ has the strong doubling property,
that is, there exists a number $\tau>1$ and a constant $A_\tau>0$ such that
for all $R>0$ and $y\in\mathbb{R}^n$ inequality \eqref{eq:strong-doubling}
holds. It is easy to see that for every open cone $C$ with the vertex at the origin, 
one can find sequences $R_j\nearrow \infty$ and 
$y_j\in C$, $j\in\mathbb{N}$ such that
\[
B(y_j,\tau R_j) \subseteq C,
\quad  
B(y_j,\tau R_j) \cap B(y_k,\tau R_k) = \emptyset 
\quad\mbox{for}\quad k \not= j. 
\]
Then \eqref{eq:strong-doubling} implies that
\[
\frac{\|\chi_{B(y_j,\tau R_j)}\|_{X(C)}}
{\|\chi_{B(y_j, R_j)}\|_{X(C)}} 
=
\frac{\|\chi_{B(y_j,\tau R_j)}\|_{X(\mathbb{R}^n)}}
{\|\chi_{B(y_j, R_j)}\|_{X(\mathbb{R}^n)}}
\le 
A_\tau , 
\quad 
j \in \mathbb{N} .
\]
So, $X(C)$ has the separated doubling property and, hence, 
the weak doubling property.
\end{proof}
%%%----------------------------------------------------------------------------
\section{Proof of Theorem~\ref{th:norm-estimate}}
\label{sec:norm-estimate-proof}
The proof is almost identical to that in \cite[Section 4.1]{KS19}. We repeat 
it here for the reader's convenience. Recall that a point $x\in\mathbb{R}^n$ 
is said to be a Lebesgue point of a function $f\in L_{\rm loc}^1(\mathbb{R}^n)$ 
if
\[
\lim_{R\to 0^+}\frac{1}{|B(x,R)|}\int_{B(x,R)}|f(y)-f(x)|\,dy=0.
\]
Given $\delta>0$ and a function $\psi$ on $\mathbb{R}^n$, we define
the function $\psi_\delta$ by
\[
\psi_\delta(\xi) := \delta^{-n} \psi(\xi/\delta),
\quad
\xi\in\mathbb{R}^n.
\]

Let $D_{X(\Omega),\varrho}$ be defined for all $\varrho>1$ by \eqref{eq:weak-doubling}.
If, for some $\varrho>1$, the quantity $D_{X(\Omega),\varrho}$ is infinite, then it is
obvious that
%%%
\begin{equation}\label{eq:norm-estimate-proof-1}
\|a\|_{L^\infty(\mathbb{R}^n)}
\le
D_{X(\Omega), \varrho} 
\|W_\Omega(a)\|_{\mathcal{B}(X_2(\Omega),X(\Omega))} .
\end{equation}
%%%
Since $X(\Omega)$ satisfies the weak doubling property, there 
exists $\varrho>1$ such that $D_{X(\Omega),\varrho}<\infty$. Take an arbitrary
Lebesgue point $\eta \in \mathbb{R}^n$ of the function $a$. Let an even
function $\varphi \in C^\infty_0(\mathbb{R}^n)$ satisfy the following
conditions:
\[
0 \le \varphi \le 1,
\quad
\varphi(x) = 1
\mbox{ for }
|x| \le 1,
\quad
\varphi(x) = 0
\mbox{ for }
|x| \ge \varrho.
\]
Let
\[
f_{\delta, \eta}(x)
:=
{e^{i\eta x}} \varphi(\delta x),
\quad
x \in \mathbb{R}^n,
\quad
\delta > 0,
\]
and
\[
f_{\delta, \eta, y}(x)
:=
f_{\delta, \eta}(x - y),
\quad
y \in \mathbb{R}^n.
\]
Then
%%%
\begin{align*}
(F f_{\delta, \eta, y})(\xi)
&=
{e^{-i\xi y}} (F f_{\delta, \eta})(\xi)
=
{e^{-i\xi y}} \delta^{-n} (F\varphi)\left(\frac{\xi - \eta}{\delta}\right)
\\
&=
{e^{-i\xi y}} (F\varphi)_\delta(\xi - \eta)
=
{e^{-i\xi y}} (F\varphi)_\delta(\eta - \xi)
\end{align*}
%%%
and
%%%
\begin{align*}
\left(F^{-1} aF f_{\delta, \eta, y}\right)(x)
&=
\frac{1}{(2\pi)^n} \int_{\mathbb{R}^n} e^{i{(x - y)}\xi} a(\xi)
(F\varphi)_\delta(\eta - \xi)\,d\xi ,
\\
a(\eta) f_{\delta, \eta, y}(x)
&=
\frac{1}{(2\pi)^n} \int_{\mathbb{R}^n}
e^{i{(x - y)}\xi} a(\eta) (F\varphi)_\delta(\eta - \xi)\,d\xi.
\end{align*}
%%%
Hence, for all $x,y\in\mathbb{R}^n$ and $\delta>0$,
%%%
\begin{align*}
&
\left|
\left(F^{-1}aFf_{\delta,\eta,y}\right)(x) - a(\eta) f_{\delta,\eta, y}(x)
\right|
\\
&\qquad=
\frac{1}{(2\pi)^n}
\left|
\int_{\mathbb{R}^n}
e^{i{(x - y)}\xi} (a(\xi) -  a(\eta)) (F\varphi)_\delta(\eta - \xi)\, d\xi
\right|
\\
&\qquad\le
\frac{1}{(2\pi)^n} \int_{\mathbb{R}^n}
|a(\xi) -  a(\eta)|\,|(F\varphi)_\delta(\eta - \xi)|\,d\xi.
\end{align*}
%%%
Since $F\varphi \in S(\mathbb{R}^n)$ and $\eta$ is a Lebesgue point of $a$,
it follows from \cite[Lemma 2.16]{KS19} that for any $\varepsilon > 0$
there exists $\delta_\varepsilon > 0$ such that for all
$x, y \in \mathbb{R}^n$ and all $\delta \in (0, \delta_\varepsilon)$,
%%%
\begin{equation}\label{eq:norm-estimate-proof-2}
\left|
\left(F^{-1} aF f_{\delta,\eta,y}\right)(x) - a(\eta) f_{\delta,\eta,y}(x)
\right| < \varepsilon.
\end{equation}
%%%
It follows from the definition of $D_{X(\Omega),\varrho}$ that there exist 
$\delta \in (0, \delta_\varepsilon)$ and $y \in \Omega$ for which
$B(y, \varrho/\delta) \subseteq \Omega$ and
%%%
\begin{equation}\label{eq:norm-estimate-proof-3}
\frac{\left\|\chi_{B(y, \varrho/\delta)}\right\|_{X(\Omega)}}
{\left\|\chi_{B(y, 1/\delta)}\right\|_{X(\Omega)}}
\le
D_{X(\Omega), \varrho} + \varepsilon.
\end{equation}
%%%
In the remaining part of the proof, we assume that $\delta$ and $y$ satisfy 
the above conditions. 

It is clear that 
%%%
\begin{equation}\label{eq:norm-estimate-proof-4}
\operatorname{supp} f_{\delta, \eta, y} 
\subseteq 
B(y, \varrho/\delta) 
\subseteq \Omega,
\quad
|f_{\delta, \eta, y}|  \chi_{B(y, 1/\delta)} = \chi_{B(y, 1/\delta)}.
\end{equation}
%%%
Then \eqref{eq:norm-estimate-proof-2} implies that 
%%%
\[
|a(\eta)| \chi_{B(y, 1/\delta)}
\le
\chi_{B(y, 1/\delta)}\left|F^{-1} aF f_{\delta,\eta,y}\right| + 
\varepsilon \chi_{B(y, 1/\delta)}   
\le |W_\Omega(a) f_{\delta,\eta,y}| + \varepsilon \chi_{B(y, 1/\delta)} .
\]
%%%
Hence
%%%
\begin{align}
&
|a(\eta)| 
\left\|\chi_{B(y, 1/\delta)} \right\|_{X(\Omega)}
\nonumber\\
 &\qquad\le
\left\|W_\Omega(a) f_{\delta,\eta,y}\right\|_{X(\Omega)}
+
\varepsilon \left\|\chi_{B(y, 1/\delta)}\right\|_{X(\Omega)}
\nonumber 
\\
&\qquad\le
\|W_\Omega(a)\|_{\mathcal{B}(X_2(\Omega),X(\Omega))}
\left\|f_{\delta, \eta, y}\right\|_{X(\Omega)}
+
\varepsilon \left\|\chi_{B(y, 1/\delta)}\right\|_{X(\Omega)}
\nonumber 
\\
&\qquad\le
\|W_\Omega(a)\|_{\mathcal{B}(X_2(\Omega),X(\Omega))}
\left\|\chi_{B(y, \varrho/\delta)}\right\|_{X(\Omega)}
+
\varepsilon \left\|\chi_{B(y, 1/\delta)}\right\|_{X(\Omega)} .
\label{eq:norm-estimate-proof-5}
\end{align}
%%%
Dividing both sides of inequality
\eqref{eq:norm-estimate-proof-5} by
$\left\|\chi_{B(y, 1/\delta)}\right\|_{X(\Omega)}$ and using 
\eqref{eq:norm-estimate-proof-3}, we get
\[
|a(\eta)| \le  (D_{X(\Omega), \varrho} + \varepsilon)
\|W_\Omega(a)\|_{\mathcal{B}(X_2(\Omega),X(\Omega))} + \varepsilon
\quad\mbox{for all}\quad \varepsilon > 0.
\]
Hence, for all Lebesgue points $\eta \in \mathbb{R}^n$ of the function $a$,
we have
\[
|a(\eta)| \le  D_{X(\Omega), \varrho} 
\|W_\Omega(a)\|_{\mathcal{B}(X_2(\Omega),X(\Omega))}.
\]
Since $a\in L^\infty(\mathbb{R}^n)\subseteq L_{\rm loc}^1(\mathbb{R}^n)$,
almost all points $\eta\in\mathbb{R}^n$ are Lebesgue points of the function
$a$ in view of the Lebesgue differentiation theorem (see, e.g.,
\cite[Theorem~1.5.1]{G24}).
Therefore, inequality
\eqref{eq:norm-estimate-proof-1} is fulfilled for all $\varrho>1$.
It is now left to apply Lemma~\ref{le:=1}(a).
\qed
%%%----------------------------------------------------------------------------
\section{Proof of Theorem~\ref{th:Kuratowski-estimate}}
\label{sec:Kuratowski-estimate-proof}
The proof is similar to the above one. In the first part, up to 
\eqref{eq:norm-estimate-proof-2}, one only needs to change 
$D_{X(\Omega), \varrho}$ by $S_{X(\Omega), \varrho}$.

It follows from the definition of $S_{X(\Omega),\varrho}$ that there exist 
$\delta_j \in (0, \delta_\varepsilon)$ and $y_j \in \Omega$, $j \in \mathbb{N}$, 
for which
\[
B(y_j, \varrho/\delta_j) \subseteq \Omega, 
\quad
B(y_j,\varrho/\delta_j) \cap B(y_k,\varrho/\delta_k) = \emptyset 
\quad\mbox{for}\quad k \not= j, 
\]
and
\begin{equation}\label{eq:Kuratowski-estimate-proof-1}
\frac{\left\|\chi_{B(y_j, \varrho/\delta_j)}\right\|_{X(\Omega)}}
{\left\|\chi_{B(y, 1/\delta_j)}\right\|_{X(\Omega)}}
\le
S_{X(\Omega), \varrho} + \varepsilon.
\end{equation}
In the remaining part of the proof, we assume that $\delta_j$ and $y_j$ 
satisfy the above conditions. 

Let 
$\varphi_{\delta_j, \eta, y_j} := 
f_{\delta_j, \eta, y_j}/\|f_{\delta_j, \eta, y_j}\|_{X(\Omega)}$. 
It follows from \eqref{eq:norm-estimate-proof-2} that
for all $x\in\mathbb{R}^n$, for all $j,k\in\mathbb{N}$, and for all
Lebesgue point $\eta$ of $a$,
%%%
\begin{align}
&
\left|
\left(F^{-1} aF (\varphi_{\delta_j,\eta,y_j} 
- 
\varphi_{\delta_k,\eta,y_k})\right)(x) 
- 
a(\eta) (\varphi_{\delta_j,\eta,y_j}(x) 
- 
\varphi_{\delta_k,\eta,y_k}(x))
\right| 
\nonumber\\
&\quad
< \varepsilon\left(
\frac{1}{\|f_{\delta_j, \eta, y_j}\|_{X(\Omega)}} 
+ 
\frac{1}{\|f_{\delta_k, \eta, y_k}\|_{X(\Omega)}}
\right). 
\label{eq:Kuratowski-estimate-proof-2}
\end{align}
%%%
It follows from \eqref{eq:norm-estimate-proof-4}
that 
\begin{equation}\label{eq:Kuratowski-estimate-proof-3}
\|\chi_{B(y_k,1/\delta_k)}\|_{X(\Omega)}
\le 
\|f_{\delta_k,\eta,y_k}\|_{X(\Omega)}
\le 
\|\chi_{B(y_k,\varrho/\delta_k)}\|_{X(\Omega)},
\quad
k\in\mathbb{N}.
\end{equation}
%%%
Since $j$ and $k$ are arbitrary natural numbers, we can assume without 
loss of generality that
%%%
\begin{equation}\label{eq:Kuratowski-estimate-proof-4}
\|\chi_{B(y_k, 1/\delta_k)}\|_{X(\Omega)} 
\le 
\|\chi_{B(y_j, 1/\delta_j)}\|_{X(\Omega)}.
\end{equation}
%%%
It follows from \eqref{eq:Kuratowski-estimate-proof-2}--%
\eqref{eq:Kuratowski-estimate-proof-4} that
%%%
\begin{align}
&
|a(\eta)|\chi_{B(y_k,1/\delta_k)}
|\varphi_{\delta_j,\eta,y_j}-\varphi_{\delta_k,\eta,y_k}|
\nonumber\\
& \qquad\le 
\chi_{B(y_k,1/\delta_k)}
|F^{-1}aF(\varphi_{\delta_j,\eta,y_j}-\varphi_{\delta_k,\eta,y_k})|
\nonumber\\
& \qquad\quad +
\varepsilon \chi_{B(y_k,1/\delta_k)}
\left(
\frac{1}{\|f_{\delta_j, \eta, y_j}\|_{X(\Omega)}} 
+ 
\frac{1}{\|f_{\delta_k, \eta, y_k}\|_{X(\Omega)}}
\right)
\nonumber\\
&\qquad \le
|W_\Omega(a)(\varphi_{\delta_j,\eta,y_j}-\varphi_{\delta_k,\eta,y_k})|
+
2\varepsilon 
\frac{\chi_{B(y_k,1/\delta_k)}}{\|\chi_{B(y_k,1/\delta_k)}\|_{X(\Omega)}}.
\label{eq:Kuratowski-estimate-proof-5}
\end{align}
%%%
It is easy to see that
%%%
\begin{equation}\label{eq:Kuratowski-estimate-proof-6}
\chi_{B(y_k, 1/\delta_k)} 
|\varphi_{\delta_j,\eta,y_j} 
- 
\varphi_{\delta_k,\eta,y_k}| 
= 
\frac{\chi_{B(y_k, 1/\delta_k)}}{\|f_{\delta_k, \eta, y_k}\|_{X(\Omega)}}.
\end{equation}
%%%
Combining \eqref{eq:Kuratowski-estimate-proof-1},
\eqref{eq:Kuratowski-estimate-proof-3},
\eqref{eq:Kuratowski-estimate-proof-5}, and
\eqref{eq:Kuratowski-estimate-proof-6}, we arrive at
%%%
\begin{align*}
&
\frac{|a(\eta)|}{S_{X(\Omega),\varrho}+\varepsilon}
\cdot
\frac{\chi_{B(y_k,1/\delta_k)}}{\|\chi_{B(y_k,1/\delta_k)}\|_{X(\Omega)}}
\\
& \qquad\le 
|a(\eta)|
\,
\frac{\chi_{B(y_k,1/\delta_k)}}{\|\chi_{B(y_k,\varrho/\delta_k)}\|_{X(\Omega)}}
\le 
|a(\eta)|
\,
\frac{\chi_{B(y_k,1/\delta_k)}}{\|f_{\delta_k,\eta,y_k}\|_{X(\Omega)}}
\\
&\qquad\le
|W_\Omega(a)(\varphi_{\delta_j,\eta,y_j}-\varphi_{\delta_k,\eta,y_k})|
+
2\varepsilon 
\frac{\chi_{B(y_k,1/\delta_k)}}{\|\chi_{B(y_k,1/\delta_k)}\|_{X(\Omega)}}.
\end{align*}
%%%
Then, by the lattice property,
\[
\frac{|a(\eta)|}{S_{X(\Omega),\varrho}+\varepsilon}
\le 
\|W_\Omega(a)(\varphi_{\delta_j,\eta,y_j}-\varphi_{\delta_k,\eta,y_k})\|_{X(\Omega)}
+
2\varepsilon. 
\]
Since $\varphi_{\delta_j,\eta,y_j}, \varphi_{\delta_k,\eta,y_k} \in B_{X_2(\Omega)}$, it follows from the above inequality that
\[
\|W_\Omega(a)\|_{\mathcal{B}(X_2(\Omega),X(\Omega)),\kappa} 
\ge 
\frac{|a(\eta)|}{2(S_{X(\Omega), \varrho} + \varepsilon)} - \varepsilon
\]
for all Lebesgue points $\eta \in \mathbb{R}^n$ of the function $a$. 
By the Lebesgue differentiation theorem (see, e.g., \cite[Theorem~1.5.1]{G24}),
the above inequality holds for a.e. $\eta\in\mathbb{R}^n$. So,
\[
\|W_\Omega(a)\|_{\mathcal{B}(X_2(\Omega),X(\Omega)),\kappa} 
\ge 
\frac{\|a\|_{L^\infty(\mathbb{R}^n)}}{2(S_{X(\Omega), \varrho} + \varepsilon)} - \varepsilon
\]
for all $\varepsilon > 0$. Hence
\[
\|W_\Omega(a)\|_{\mathcal{B}(X_2(\Omega),X(\Omega)),\kappa} 
\ge 
\frac{\|a\|_{L^\infty(\mathbb{R}^n)}}{2S_{X(\Omega), \varrho}}.
\]
It is now left to apply Lemma~\ref{le:=1}(b).
\qed
%%%----------------------------------------------------------------------------
\section{The case of weighted variable Lebesgue spaces}
\label{sec:weighted-variable-Lebesgue-spaces}
It is instructive to look at Corollary~\ref{co:nice-for-cones} in the case 
of weighted variable Lebesgue spaces. Let $p(\cdot):\mathbb{R}^n\to(1,\infty)$ 
be a measurable function, called a variable exponent. For simplicity, we 
assume that
%%% 
\begin{equation}\label{eq:non-trivial-exponents}
1<
p_-:=\operatornamewithlimits{ess\,inf}_{x\in\mathbb{R}^n} p(x),
\quad
p_+:=\operatornamewithlimits{ess\,inf}_{x\in\mathbb{R}^n} p(x)
<\infty.
\end{equation}
%%%
In this case, the conjugate exponent $p'(\cdot):\mathbb{R}^n\to(1,\infty)$
is well defined by
\[
1/p(x)+1/p'(x)=1,
\quad
x\in\mathbb{R}^n.
\]
The variable Lebesgue space $L^{p(\cdot)}(\mathbb{R}^n)$
consists of all complex-valued measurable functions $f$ on $\mathbb{R}^n$
such that
\[
m_{p(\cdot)}(f/\lambda):=\int_{\mathbb{R}^n} |f(x)/\lambda|^{p(x)} dx<\infty
\]
for some $\lambda=\lambda(f)>0$. It is well known that 
$L^{p(\cdot)}(\mathbb{R}^n)$ is a Banach function space
equipped with the norm
\[
\|f\|_{L^{p(\cdot)}(\mathbb{R}^n)}
:=
\inf\left\{\lambda>0\ : \ m_{p(\cdot)}(f/\lambda)\le 1\right\}.
\]
and the associate space of $L^{p(\cdot)}(\mathbb{R}^n)$ is the space 
$L^{p'(\cdot)}(\mathbb{R}^n)$ (up to equivalence of the norms)
(see \cite[Theorems~2.17, 2.32 and Section~2.10.3]{CF13}
or \cite[Theorems~3.2.7 and~3.2.13]{DHHR11}).
The hypothesis \eqref{eq:non-trivial-exponents} ensures that 
$L^{p(\cdot)}(\mathbb{R}^n,w)$ is separable and reflexive
(see \cite[Theorem~2.78, Corollary~2.81]{CF13} or 
\cite[Corollary~3.4.5 and Theorem~3.47]{DHHR11}).

In the case of weighted variable
Lebesgue spaces, Corollary~\ref{co:nice-for-cones} reads as follows.
%%%----------------------------------------------------------------------------
\begin{theorem}\label{th:case-of-weighted-variable-Lebesgue-space}
Let $p(\cdot):\mathbb{R}^n\to(1,\infty)$ be a variable exponent satisfying
\eqref{eq:non-trivial-exponents} and $w$ be a weight satisfying
%%%
\begin{equation}\label{eq:Muckenhoupt}
\sup_{B}\frac{1}{|B|}
\|w\chi_B\|_{L^{p(\cdot)}(\mathbb{R}^n)}
\|w^{-1}\chi_B\|_{L^{p'(\cdot)}(\mathbb{R}^n)}
<\infty,
\end{equation}
%%%
where the supremum is taken over all balls $B\subseteq\mathbb{R}^n$.
Suppose $C\subseteq\mathbb{R}^n$ is an open cone with the vertex at the origin.
If $a \in \mathcal{M}^0_{L^{p(\cdot)}(\mathbb{R}^n,w)}$, then
\[
\|a\|_{L^\infty(\mathbb{R}^n)} 
\le 
\|W_C(a)\|_{\mathcal{B}(L^{p(\cdot)}(C,w))},
\quad
\frac12\, \|a\|_{L^\infty(\mathbb{R}^n)} 
\le 
\|W_C(a)\|_{\mathcal{B}(L^{p(\cdot)}(C),w),\kappa} .
\]
In particular,
%%%
\begin{equation}\label{eq:noncompacntess-equivalence-1}
W_{\mathbb{R}^n}(a)\in\mathcal{K}(L^{p(\cdot)}(\mathbb{R}^n,w))
\ \Longleftrightarrow\
a=0
\end{equation}
%%%
and
\begin{equation}\label{eq:noncompactness-equivalence-2}
W_{C}(a)\in\mathcal{K}(L^{p(\cdot)}(C,w))
\ \Longleftrightarrow\ 
a=0.
\end{equation}
\end{theorem}
%%%----------------------------------------------------------------------------
\begin{proof}
If \eqref{eq:Muckenhoupt} holds, then 
$w\in L_{\rm loc}^{p(\cdot)}(\mathbb{R}^n)$
and $w^{-1}\in L_{\rm loc}^{p'(\cdot)}(\mathbb{R}^n)$.
Since $L^{p(\cdot)}(\mathbb{R}^n)$ is a separable and reflexive
Banach function space, it follows from Lemma~\ref{le:wBFS} that
\[
L^{p(\cdot)}(\mathbb{R}^n,w)=
\left\{
f\in\mathfrak{M}(\mathbb{R}^n)\ :\ fw\in L^{p(\cdot)}(\mathbb{R}^n)
\right\}
\]
is a separable and reflexive Banach function space equipped with the 
norm $\|f\|_{L^{p(\cdot)}(\mathbb{R}^n)}=\|fw\|_{L^{p(\cdot)}(\mathbb{R}^n)}$,
and its associate space $(L^{p(\cdot)})'(\mathbb{R}^n,w)$
coincides with $L^{p'(\cdot)}(\mathbb{R}^n,w^{-1})$
(up to equivalence of the norms). Hence, condition \eqref{eq:Muckenhoupt} is 
equivalent to condition \eqref{eq:Berezhnoi}. It remains to apply
Corollary~\ref{co:nice-for-cones}.
\end{proof}
%%%----------------------------------------------------------------------------
It is easy to see that if $p(\cdot):\mathbb{R}^n\to(1,\infty)$ is the constant
equal to $p$, then $L^{p(\cdot)}(\mathbb{R}^n,w)$ is nothing but the standard
weighted Lebesgue space $L^p(\mathbb{R}^n,w)$ and \eqref{eq:Muckenhoupt} is 
nothing but the classical Muckenhoupt condition: 
\[
\sup_{B}
\left(\frac{1}{|B|}\int_B w^p(x)\,dx\right)^{1/p}
\left(\frac{1}{|B|}\int_B w^{-p'}(x)\,dx\right)^{1/p'},
\]
where the supremum is taken over all balls $B\subseteq\mathbb{R}^n$
(see, e.g., \cite[Ch.~8]{G24}).
The equivalence in 
\eqref{eq:noncompacntess-equivalence-1} was established for $n=1$ 
and $p(\cdot)=\operatorname{const}=p$ in \cite[Corollary~1.2]{FKK19}.
We were not able to find the equivalence in 
\eqref{eq:noncompactness-equivalence-2} in this case for 
$C=\mathbb{R}_+$ explicitly stated in the literature.

If $n=1$ and $w=1$, then it follows from 
\cite[Theorems~5.3.4 and 5.7.2]{DHHR11} that condition \eqref{eq:Muckenhoupt}
is not sufficient for the boundedness of the Hardy-Littlewood maximal
operator on a variable Lebesgue space $L^{p(\cdot)}(\mathbb{R})$ without
additional assumptions on the variable exponent
(see \cite[Theorem~1.3]{CDH11} and \cite[Theorem~1.5]{CFN12} for conditions on
$p(\cdot):\mathbb{R}^n\to(1,\infty)$ guaranteeing that \eqref{eq:Muckenhoupt} 
is equivalent to the boundedness of $M$ on $L^{p(\cdot)}(\mathbb{R}^n,w)$).
So, in the case $n=1$ and $w=1$, the equivalence in 
\eqref{eq:noncompacntess-equivalence-1} does not follow from 
\cite[Theorem~1.1]{FKK19}.

Let us conclude the paper by stating the corollary of our results under the 
natural assumptions used in the study of Wiener-Hopf operators on
variable Lebesgue spaces (see \cite{KV24,V-JMS}), which gives 
a satisfactory answer to Question~\ref{question:main}.
%%%----------------------------------------------------------------------------
\begin{corollary}
Let $p(\cdot):\mathbb{R}^n\to(1,\infty)$ be a variable exponent satisfying
\eqref{eq:non-trivial-exponents} and such that the Hardy-Littlewood maximal
operator $M$ is bounded on the variable Lebesgue space 
$L^{p(\cdot)}(\mathbb{R}^n)$. If
$a\in\mathcal{M}_{L^{p(\cdot)}(\mathbb{R}^n)}^0$ and $C\subseteq\mathbb{R}^n$
is an open cone with the vertex at the origin, then $W_C(a)$ is compact on
$L^{p(\cdot)}(C)$ if and only if $a=0$.
\end{corollary}
%%%----------------------------------------------------------------------------
Indeed, if $M$ is bounded on $L^{p(\cdot)}(\mathbb{R}^n)$, then it follows 
from \cite[Theorem~4.4.10, the discussion before it, and Lemma~4.4.13]{DHHR11}
that condition \eqref{eq:Muckenhoupt} is fulfilled with $w=1$.
So, the desired result is a consequence of 
Theorem~\ref{th:case-of-weighted-variable-Lebesgue-space}.
%%%----------------------------------------------------------------------------
\subsection*{Funding}
This work is funded by national funds through the FCT – Fundação para a 
Ciência e a Tecnologia, I.P., under the scope of the projects 
UID/00297/2025 (\url{https://doi.org/10.54499/UID/00297/2025}) and 
UID/PRR/00297/2025 
(\url{https://doi.org/10.54499/UID/PRR/00297/2025}) .
%%%----------------------------------------------------------------------------
\bibliographystyle{abbrv}
\bibliography{OKES26-arXiv}
\end{document}